\documentclass{article}

\usepackage{authblk}
\usepackage[english]{babel}
\usepackage{amssymb, amsmath}
\usepackage{wrapfig}
\usepackage{varioref}
\usepackage{graphicx}
\usepackage{caption}
\usepackage{subcaption}
\usepackage{adjustbox}
\usepackage{placeins}
\usepackage{url}
\usepackage[utf8]{inputenc}
\usepackage[english]{babel}
\usepackage[normalem]{ulem}
\usepackage[nolist]{acronym}
\usepackage[dvipsnames]{xcolor}
\usepackage{colortbl}
\definecolor{GrayLight}{gray}{0.9}
\definecolor{GrayDark}{gray}{0.75}
\usepackage{lipsum}
\usepackage{titlesec}
\usepackage{dsfont}
\usepackage{multirow}
\usepackage{svg}
\usepackage{url}
\usepackage{hyperref}
\hypersetup{
  colorlinks=true,
  linkcolor=blue,
}
\usepackage{tikz}
\usepackage{supertabular}
\usepackage{booktabs}
\usetikzlibrary{decorations.pathreplacing}
\addto\captionsenglish{}
\usepackage{framed}
\usepackage{scrextend}
\addtokomafont{labelinglabel}{\sffamily\bfseries}

\newcommand{\m}[1]{\ensuremath{#1}}

\newcommand{\ns}[2]{_{#1}^{\mathrm{#2}}}

\newcommand{\Time}{\m{\mathcal{T}}}
\newcommand{\Sce}{\m{\mathcal{S}}}
\newcommand{\forallT}{\m{\forall \, t \in \mathcal{T}}}
\newcommand{\forallTS}{\m{\forall \, t \in \mathcal{T}, \, s \in \mathcal{S}}}

\newcommand{\Demand}[2]{\m{D\ns{#1}{#2}}}
\newcommand{\PV}[2]{\m{PV\ns{#1}{#2}}}
\newcommand{\PVpeak}{\m{PV^{peak}}}
\newcommand{\Bat}[2]{\m{B\ns{#1}{#2}}}
\newcommand{\ElecPrice}[2]{\m{\lambda\ns{#1}{#2}}}
\newcommand{\Sto}[2]{\m{S\ns{#1}{#2}}}
\newcommand{\Prob}[1]{\m{\pi_{#1}}}

\newcommand{\bat}[2]{\m{b\ns{#1}{#2}}}
\newcommand{\powerTransac}[2]{\m{x\ns{#1}{#2}}}
\newcommand{\sto}[2]{\m{s\ns{#1}{#2}}}
\newcommand{\pv}[2]{\m{pv\ns{#1}{#2}}}
\newcommand{\grid}[2]{\m{g\ns{#1}{#2}}}

\newcommand{\deltaT}{\m{\Delta t}}

\newcommand{\toAC}{\m{\phi^{\text{DC/AC}}}}
\newcommand{\toDC}{\m{\phi^{\text{AC/DC}}}}
\newcommand{\charge}{\m{\eta^{+}}}
\newcommand{\discharge}{\m{\eta^{-}}}

\title{Economic evaluation of stochastic home energy management systems in a
  realistic rolling horizon setting}

\author[*,1,3]{ Julian Lemos-Vinasco}

\author[1]{ Amos Schledorn}
\author[2]{ S. Ali Pourmousavi}
\author[1]{ Daniela Guericke}

\affil[*]{\small Corresponding author, 
    \textit{Email address:} jlvi@dtu.dk, julian.lemos@watts.dk,
    \textit{Postal Address:} Anker Engelunds Vej 1, Building 101A, 2800 Kongens
    Lyngby, Denmark}
\affil[1]{\small Technical University of Denmark (DTU), Denmark}
\affil[2]{\small The University of Adelaide, Australia}
\affil[3]{\small Watts A/S, Denmark}

\begin{acronym}
  \acro{DA}{day-ahead}
  \acro{AC}{Alternate Current}
  \acro{DC}{Direct Current}
  \acro{PV}{photovoltaic}
  \acro{RES}{renewable energy source}
  \acro{TS}{time series}
  \acro{EL}{electrical load}
  \acro{HEMS}{home energy management system}
  \acro{SoC}{State of Charge}
  \acro{LF}{load forecasting}
  \acro{PLF}{probabilistic load forecasting}
  \acro{NWP}{numerical weather prediction}
  \acro{SDE}{Stochastic Differential Equation}
  \acro{RLS}{Recursive Least Squares}
  \acro{MILP}{mixed-integer linear program}
  \acro{AR}{autoregressive}
  \acro{SHT}{smart home technologies}
\end{acronym}

\begin{document}
\maketitle
\begin{abstract}
  Home energy management systems (HEMSs) are expected to become a crucial part
  of future smart grids. However, there is a limited number of studies that
  comprehensively assess the potential economic benefits of \acp{HEMS} for
  consumers under real market conditions and which take account of consumers'
  capabilities. In this study, a new optimization-based \ac{HEMS} controller is
  presented to operate a photovoltaic and battery system. The \ac{HEMS}
  controller considers the consumers' electrical load uncertainty by integrating
  multivariate probabilistic forecasting methods and a stochastic optimization
  in a rolling horizon. As a case study, a comprehensive simulation study is
  designed to emulate the operation of a real \ac{HEMS} using real data from
  nine Danish homes over different seasons under real-time retail prices. The
  optimization-based control strategies are compared with a default (naive)
  control strategy that encourages self consumption. Simulation results show
  that seasonality in the consumers' load and electricity prices have a
  significant impact on the performance of the control strategies. A combination
  of optimization-based and naive control strategy presents the best overall
  results.

\vspace{5pt}
\end{abstract}

\pagebreak
\subsection*{Nomenclature}
\tablefirsthead{\toprule}
\tablehead{\midrule}
\tabletail{\midrule}
\tablelasttail{\bottomrule}

\footnotesize
\begin{supertabular}{p{0.14\linewidth}p{0.8\linewidth}}
 
  \multicolumn{2}{l}{Sets}\\\midrule
  $\Time$ & Set of time steps $t$\\
  $\Sce$ & Set of scenarios $s$\\

  \midrule\multicolumn{2}{l}{Parameters} \\\midrule
  $\Demand{s,t}{}$ & Electricity load in scenario $s \in \Sce$ and period
  $t \in \Time$ [kWh]
  \\
  $\ElecPrice{t}{+}$ & Electricity purchase cost in period $t \in \Time$
  [DKK/kWh]
  \\
  $\ElecPrice{t}{-}$ & Electricity sale price in period $t \in \Time$
  [DKK/kWh]
  \\
  $\PV{t}{}$ & \ac{PV} production in period $t \in \Time$ [kWh]
  \\
  $\PVpeak$ & \ac{PV} system peak production [Wh]
  \\
  $\Sto{}{ini}$ & Initial battery \ac{SoC} [Kwh]
  \\
  $\Sto{}{max}$ & Battery maximum storage capacity [kWh]
  \\
  $\Sto{}{min}$ & Battery minimum storage capacity [kWh]
  \\
  $\Bat{}{in}$ & Battery charge limit per period [kW]
  \\
  $\Bat{}{out}$ & Battery discharge limit period [kW]
  \\
  $\Prob{s}$ & Probability of scenario $s \in \Sce$
  \\
  $\toAC$ & Efficiency factor when inverting power flows from \ac{DC} to
  \ac{AC}
  \\
  $\toDC$ & Efficiency factor when converting power flows from \ac{AC} to
  \ac{DC}
  \\
  $\charge$ & Battery charge efficiency factor
  \\
  $\discharge$ & Battery discharge efficiency factor
  \\
  $M$ & Big $M$ value define as $M = \max{(\PV{t}{})} + \max{(\Demand{s,t}{})} + \Bat{}{in}$
  \\
  \midrule\multicolumn{2}{l}{Variables} \\\midrule
  $\powerTransac{s,t}{+} \in \mathbb{R}^{+}$ & Electricity bought from the
  electricity retailer in scenario $s \in \Sce$ and period $t \in \Time$ [kWh]
  \\
  $\powerTransac{s,t}{-} \in \mathbb{R}^{+}$ & Electricity  sold to the
  electricity retailer in scenario $s \in \Sce$ and period $t \in \Time$ [kWh]
  \\
  $\grid{s,t}{D} \in \mathbb{R}^{+}$ & Power from the grid used to satisfy the demand  in scenario $s \in \Sce$ and period $t \in \Time$ [kW]
  \\
  $\grid{s,t}{b} \in \mathbb{R}^{+}$ & Power sent from the grid to the battery
  in scenario $s \in \Sce$ and period $t \in \Time$ [kW]
  \\
  $\bat{s,t}{+} \in \mathbb{R}^{+}$ & Battery charge power flow in scenario
  $s \in \Sce$ and period $t \in \Time$ [kW]
  \\
  $\bat{s,t}{-} \in \mathbb{R}^{+}$ & Battery discharge power flow in scenario
  $s \in \Sce$ and period $t \in \Time$ [kW]
  \\
  $\bat{t,s}{g} \in \mathbb{R}{+}$ & Power delivered from the battery to the grid
   in scenario $s \in \Sce$ and period $t \in \Time$ [kW]
  \\
  $\bat{t,s}{D} \in \mathbb{R}{+}$ & Power from the battery  used to satisfy the
  demand in scenario $s \in \Sce$ and period $t \in \Time$
  [kW]
  \\
  $\pv{t}{g} \in \mathbb{R}{+}$ & Power delivered directly from the \ac{PV} system to the grid in period $t \in \Time$ [kW]
  \\
  $\pv{t}{D} \in \mathbb{R}{+}$ & Power from the \ac{PV} system used to satisfy the demand in period $t \in \Time$ [kW]
  \\
  $\pv{t}{b} \in \mathbb{R}{+}$ & Power from the \ac{PV} system  to the battery in period $t \in \Time$ [kW]
  \\
  $s_{s,t} \in \mathbb{R}^{+}$ & Battery \ac{SoC} in scenario $s \in \Sce$ and
  period $t \in \Time$ [kWh]
  \\
  $z_{s,t} \in \{0,1\}$ & Binary variable, indicating if electricity was purchased or sold
  to the grid for scenario $s \in \Sce$ and period $t \in \Time$
  \\
  $y_{s,t} \in \{0,1\}$ & Binary variable, indicating if the battery is charging or
  discharging for scenario $s \in \Sce$ and period $t \in \Time$
  \\
\end{supertabular}
\normalsize

\section{Introduction}
As one of the major smart grid technologies, home energy management systems
(\acp{HEMS}) are expected to play a key role managing energy consumption at the
residential level by reacting to real-time prices and/or $\text{CO}_{2}$-based
signals. In Europe, this coincides with efforts of electricity market operators
and policy makers to push for a wider adoption of real-time tariffs for
residential consumers that reflect the true condition of the power system and
provide cost-savings for consumers \cite{Europeanparliament,Energinet}.

High expectations have been placed on \acp{HEMS} by many industry stakeholders
given the systems' potential to provide a dynamic combination of production,
storage, and flexible demand \cite{Mariano-Hernandez2021,Goulden2014,Hansson}.
Therefore, studies on this topic have emerged from a variety of disciplines over
the last decade, focusing on different components of the \acp{HEMS}. Typically,
\acp{HEMS} rely on a combination of \ac{SHT} such as smart meters, sensing
devices, communication hardware and protocols, smart appliances, controllers,
and optimization techniques \cite{Balakrishnan2021}. The operation and
coordination of these components entails technical difficulties, especially for
\ac{SHT} that depend on manual intervention from end users. This has led to a
literature bias towards \ac{SHT} solutions that require minimal consumer
intervention \cite{McIlvennie2020a}.

In this regard, several studies have proposed sophisticated technical solutions
by assuming a direct control of several \ac{SHT}. These solutions have been used
for direct control of the heating systems of homes and buildings, smart
appliances, \acp{RES}, batteries, and electrical vehicle chargers (in both
grid-2-vehicle and vehicle-2-grid modes), with some parameters being defined by
the consumer \cite{Yousefi2020,Salgado2021}. Furthermore, it is important for
\acp{HEMS} to consider complex system features such as the multi-seasonality,
non-stationarity, and stochasticity of \acp{RES} and consumers' \ac{EL}
\cite{Lemos-Vinasco2021a,Ebrahimi2022}. Therefore, recent studies have included
several of these features. In \cite{Shafie-Khah2018}, a stochastic \ac{HEMS} was
proposed that considered consumers' satisfaction cost and fatigue towards
demand-response signals. The authors of \cite{Zeynali2020} proposed a two-stage
stochastic model with scenarios for wind power and electric vehicles'
availability. In \cite{Correa-Florez2018}, a similar approach is used with
additional considerations for the battery degradation cost. Other studies apply
rolling horizon approaches. Such approaches provide an opportunity to
re-optimize the problem when new information about stochastic elements are
available, for example, \ac{PV} forecast \cite{Paterakis2016a,Elkazaz2019a}.

The studies in the literature on control strategies for \acp{HEMS} display
several similarities. First, the studies assume direct control over different
\ac{SHT} (controlled laboratory conditions and/or simulations). Second, most
studies are mainly oriented towards demand-response programs by assuming access
to the wholesale market electricity prices (day-ahead and/or intra-day prices).
Third, they present a cost-benefit comparison for a limited time period (ranging
from days to weeks), typically in cold seasons with a passive consumer (consumer
without \ac{SHT}) as the baseline. In contrast to these publications, the
results of field studies and trial projects have questioned the real benefits
that consumers will be able to perceive. Results from a nine-month field trial
with ten households in the UK concluded that ``there is little evidence that
\ac{SHT} will generate substantial energy saving and, indeed, there is a risk
that they may generate a form of energy intensification'' \cite{Hargreaves2017}.
These observations are aligned with the findings in \cite{Nicholls2017}, where
in a trial with 40 households with basic \ac{SHT}, minimum economic benefits
were reported, with some households reporting energy intensification. However,
the setups of these studies used smart appliances requiring manual consumer
interventions. Moreover, \cite{Nicholls2015} suggest that very limited economic
benefits can be expected from these types of setups because of the inherent
inflexibility of some consumers.

Although the above studies indicate a need for more research on \ac{SHT} that
require manual intervention from consumers, the main body of literature assumes
direct access and control of most \ac{SHT} elements. This is a strong assumption
that may distort the studies' results \cite{TiradoHerrero2018}. Furthermore, the
results mainly describe technical aspects with assumptions that may not work
under current market rules. For instance, in
\cite{Correa-Florez2018,Paterakis2016a,Yousefi2020}, it is not clear if the
electricity prices used correspond to prices accessible to consumers, or they
assumed that consumers have access to the wholesale electricity markets, which
is not possible due to the small size of individual consumers' load and
\acp{RES} in the European markets \cite{Nordpool_regulations}. Assuming access
to wholesale market prices disregards the fact that consumers are subject to
taxes, levies, and fees, which may have a significant impact on the results.
Additionally, the cost comparisons are made with a passive consumer as baseline,
disregarding the fact that simple self-consumption control strategies have
proved to bring significant cost reductions \cite{Azuatalam2019a}.

On the basis of the above discussion, one can argue that there is a research gap
in relation to the assessment of the economic potential of \ac{HEMS} under
system conditions and market rules accessible to the residential consumers.
These conditions must include a realistic \ac{HEMS} setup, end-consumer prices,
cover a substantial period of time consisting of different seasons, and compare
the results to self-consumption control strategies.

Based on the research gap, this paper contributes with a comprehensive economic
assessment of a \ac{HEMS} under realistic consumer and electricity market
conditions over different seasons. We propose a novel \ac{HEMS} control strategy
that uses stochastic optimization framework in a rolling horizon approach and
probabilistic forecast. At a technical level, this paper also contributes by
integrating two different multivariate probabilistic forecast methods which
consider temporal correlation and a \ac{HEMS} setup modeled as a stochastic
\ac{MILP}. Moreover, the rolling horizon approach is used to allow the
possibility of re-optimizing according to the latest information available to
the system. The data used in the case study corresponds to nine households
located in Copenhagen, Denmark, together with real hourly electricity retail
prices offered by a utility company. Furthermore, a \ac{HEMS} setup with only a
\ac{PV} and battery system is considered to emulate the possibilities that most
residential consumers have at present. Although electric vehicles are a key
element of the \acp{HEMS} of the future, the adoption of electric vehicles is
still low in Denmark \cite{Sevdari2021} and therefore they were not included in
the analysis. Operational and cost results of the proposed optimization-based
strategies are compared with a passive consumer as well as a self-consumption
(naive) control strategy.

Overall, the results indicate that a combination of an optimization-based and a
naive control strategy presents a higher economic benefit for residential
consumers throughout the year. Key research findings are summarized below:
\begin{enumerate}
  \item The stochasticity of the consumers' \ac{EL} has a significant impact on
        the performance of the optimization-based control strategies.
  \item Strong seasonality in consumption patterns shows a significant effect on
        the assessment and selection of the control strategies, with a
        self-consumption strategy (naive control) outperforming the
        optimization-based controllers in spring and summer.
  \item Under current market rules, residential consumers are not sufficiently
        incentivized to actively participate in the electricity market besides
        covering their electricity demand.
\end{enumerate}

This paper starts by presenting the \ac{HEMS} setup and the mathematical details
of the implemented models in Section \ref{sec:methods}. Next, the data and the
case study are explained in Section \ref{sec:case_study}. The simulation results
are presented in Section \ref{sec:results}, which includes a comparison between
different control strategies, and a comprehensive cost analysis. Finally, a
discussion of the findings and perspectives for future work are outlined in
Section \ref{sec:discussion}. The paper is concluded in Section
\ref{sec:conclusion}.

\section{Modeling and optimization of \ac{HEMS}}
\label{sec:methods}
\subsection{\ac{HEMS} setup}
\label{sec:hems_setup}

\begin{figure*}[t]
  \centering
  \includegraphics[width=\textwidth]{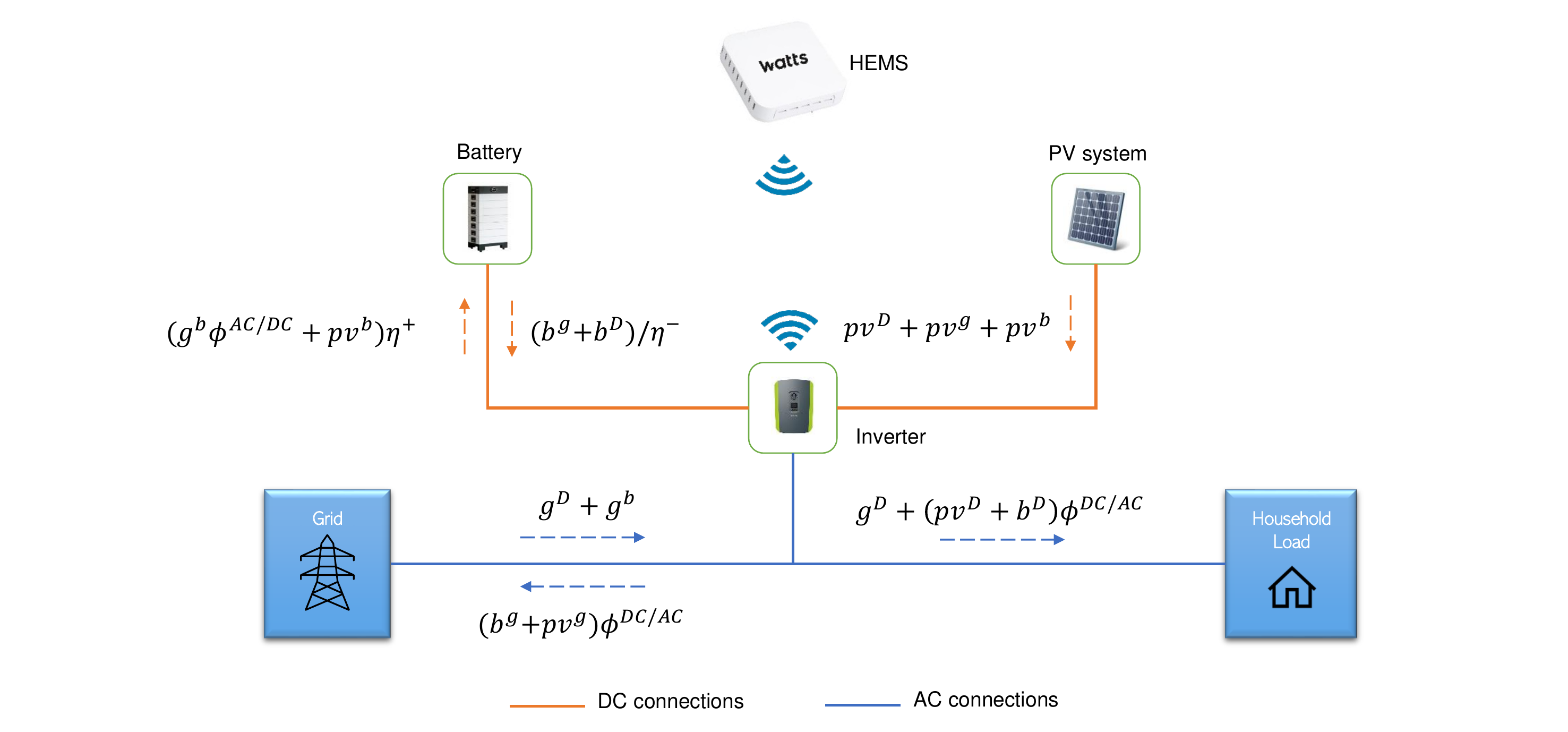}
  \caption{Schematic diagram of the home's setup showing the \ac{AC} and
    \ac{DC} connections. The different average power flow terms, which are used in the
    mathematical representation of the system (Section \ref{sec:hems_math}), are
    also included.}
  \label{fig:HEMS_schematic}
\end{figure*}

The \ac{HEMS} setup considered in this paper reflects the current network
conditions in Denmark that residential electricity consumers with access to a
\ac{PV} and a home battery system face. Moreover, a minimal control approach is
considered for the \ac{HEMS} model. This means that the \ac{HEMS} has direct
control of the home battery, but it does not have direct control over the home
appliances. A graphical overview of the setup is given in Figure
\ref{fig:HEMS_schematic}.

The electricity generation from the \ac{PV} system can be used to charge the
home battery, to meet the \ac{EL} demand, or can be exported to the grid. The
electricity losses due to AC/DC and DC/AC conversions are included in the
formulation. We assume that the electricity retailer communicates price
information to the \ac{HEMS} and that the \ac{HEMS} has access to \acp{NWP}.
Furthermore, \acp{NWP} are input to \ac{PLF} models used for the creation of
\ac{EL} scenarios. The data of real consumers in Denmark is used, however, these
consumers did not have \ac{PV} installations. Thus, a simulation model for
\ac{PV} production is implemented and presented in Section~\ref{sec:pv_model}.
The remainder of this section introduces the mathematical model formulation for
the above setup.

\subsection{HEMS optimization model}
\label{sec:hems_math}

The \ac{HEMS} controller is formulated as a stochastic \ac{MILP}
\cite{Birge2011}, where the \ac{EL} is the only uncertain parameter, i.e. having
varying realizations across scenarios. The \ac{MILP} in \eqref{eq:mip} minimizes
the expected cost for fulfilling the \ac{EL} in all scenarios $\Sce$ and periods
$\Time$. The main decision variables are the flows between the \ac{PV}, grid,
and battery components. The model considers several time periods due to the
temporal interdependence imposed by the battery \ac{SoC} in a rolling horizon
manner. This means that, when applying the solution to the \ac{HEMS}, only the
optimal solution for the first time period is applied in practice. Decisions in
subsequent periods are only considered to find optimal decisions for the first
time step. This allows for re-optimization and taking relevant decisions with
updated forecasts.
{\allowdisplaybreaks
\small
\begin{subequations}
\begin{gather}
    \label{eq:obj_func}
    \min_{\powerTransac{s,t}{-}, \powerTransac{s,t}{+}}\quad \sum_{s \in \Sce} \sum_{t \in \Time}
    \Prob{s} (\ElecPrice{t}{+} \powerTransac{s,t}{+} - \ElecPrice{t}{-} \powerTransac{s,t}{-})
\end{gather}
subject to
\begin{align}
  \label{eq:satisfy_demand}
  &\Demand{s,t}{} = (\grid{s,t}{D} + (\pv{t}{D} + \bat{s,t}{D}) \toAC) \deltaT &\forallTS
  \\
  \label{eq:purchase_balance}
  &\powerTransac{s,t}{+} = (\grid{s,t}{D} + \grid{s,t}{b}) \deltaT &\forallTS
  \\
  \label{eq:sale_balance}
  &\powerTransac{s,t}{-} = (\bat{s,t}{g} + \pv{t}{g}) \toAC \deltaT &\forallTS
  \\
  \label{eq:binary1}
  &\powerTransac{s,t}{-} \leq M z_{s,t} &\forallTS
  \\
  \label{eq:binary2}
  &\powerTransac{s,t}{+} \leq M (1-z_{s,t}) &\forallTS
  \\
  \label{eq:pv_balance}
  &\PV{t}{} = (\pv{t}{g} + \pv{t}{D} + \pv{t}{b}) \deltaT &\forallT
  \\
  \label{eq:battery_charge}
  &\bat{s,t}{+} = (\grid{s,t}{b} \toDC + \pv{t}{b})\charge &\forallTS
  \\
  \label{eq:battery_discharge}
  & \bat{s,t}{-} = (\bat{s,t}{g} + \bat{s,t}{D})/\discharge &\forallTS
  \\
  \label{eq:power_capacity}
  & \bat{s,t}{+} \leq \Bat{}{in} &\forallTS
  \\
  \label{eq:power_capacity_1}
  & \bat{s,t}{-} \leq \Bat{}{out} & \forallTS
  \\
  \label{eq:battery_binary_1}
  & \bat{s,t}{+} \leq M y_{s,t} &\forallTS
  \\
  \label{eq:battery_binary_2}
  & \bat{s,t}{-} \leq M (1 - y_{s,t}) &\forallTS
  \\
  \label{eq:storage1}
  &\sto{s,t}{} = \Sto{}{ini} + (\bat{s,t}{+} - \bat{s,t}{-}) \deltaT &\forall \, s \in \Sce, \, t=1
  \\
  \label{eq:storage2}
  &\sto{s,t}{} = \sto{s,t-1}{} + (\bat{s,t}{+} - \bat{s,t}{-}) \deltaT &\forall \, t \in \Time, \, s \in  \Sce, \, t > 1
  \\
  \label{eq:storage_capacity}
  &\Sto{}{min} \leq \sto{s,t}{} \leq \Sto{}{max} &\forallTS
  \\
  \label{eq:first_stage1}
  &\bat{s,t}{+} = \bat{j,t}{+}\, , \, \bat{s,t}{-} = \bat{j,t}{-} &\forall  s, j \in \Sce, \,  s \ne j, \, t=1
  \\
  \label{eq:first_stage2}
  &\sto{s,t}{} = \sto{j,t}{} &\forall s, j \in \Sce, \, s \ne j, \, t=1
\end{align}
\label{eq:mip}
\end{subequations}}
\normalsize

\noindent The objective function \eqref{eq:obj_func} minimizes the expected cost
of electricity. Furthermore, constraint \eqref{eq:satisfy_demand} ensures that
the consumer's demand is satisfied in all scenarios. Electricity purchase and
sale quantities are set in constraints \eqref{eq:purchase_balance} and
\eqref{eq:sale_balance}. Constraints \eqref{eq:binary1} and \eqref{eq:binary2}
ensure that electricity sale and purchase are mutually exclusive. The PV power
balance is set in constraint \eqref{eq:pv_balance} such that the total
generation meets the sum of \ac{PV} production to grid, demand and battery.
Constraints \eqref{eq:battery_charge}, \eqref{eq:battery_discharge},
\eqref{eq:power_capacity}, and \eqref{eq:power_capacity_1} model the physical
battery behaviour in terms of power flow. Simultaneous charging and discharging
of the battery is disallowed in constraints \eqref{eq:battery_binary_1} and
\eqref{eq:battery_binary_2}. The evolving \ac{SoC} is modelled by constraints
\eqref{eq:storage1} and \eqref{eq:storage2}, while constraint
\eqref{eq:storage_capacity} limits the \ac{SoC} to the battery capacity.

Since it is possible to re-optimize the solution after one time period, and the
energy exchange with the grid is unrestricted, we can frame the problem as a
two-stage stochastic problem. The first stage of the problem defines the
operational schedules of the battery in the first period, as given by
constraints \eqref{eq:first_stage1} and \eqref{eq:first_stage2}. Thus, the
battery charging and discharging in the first time period needs to be the same
for all scenarios.

\subsection{Electrical load forecast}
\label{sec:load_forecast_description}

The \ac{HEMS} optimization model presented in Section~\ref{sec:hems_math} uses
\ac{EL} scenarios as input. The scenarios must consider the temporal correlation
inherent to the \ac{EL}. Thus, the multivariate \ac{PLF} methods presented in
\cite{Lemos-Vinasco2021a} are used in this paper to generate the required
scenarios. The methods use either \ac{RLS} with a full covariance model for the
residuals or the quantile-copula with a full covariance model of the temporal
correlation under the Gaussian domain and are referenced in
\cite{Lemos-Vinasco2021a} as \emph{RLS-Free} and \emph{Copula-Free}.

\subsection{\ac{PV} simulation}
\label{sec:pv_model}

Another key element of the \ac{HEMS} is the \ac{PV} system. In this study,
\ac{PV} generation data were not available. Thus, a simulation approach is used
to estimate a rooftop \ac{PV} production. The simulation model is based on the
guidelines provided in the energy data catalogue by the Danish Energy Agency
\cite{DanishEnergyAgency}. The report suggests that electricity produced by a
\ac{PV} system should be estimated as \small
\begin{align}
  \PV{t}{} = \PVpeak \cdot \frac{GHI_{t}}{1000} \cdot \PV{}{tf} \cdot \PV{}{pr}
\end{align}
\normalsize

\noindent where \PVpeak corresponds to the \ac{PV} production under laboratory
standard test conditions (1000 $W/m^{2}$ irradiation with a cell temperature of
25\textdegree{C}), $PV^{tf}$ is the \ac{PV} transposition factor and $PV^{pr}$
is the \ac{PV} performance ratio. Moreover, the $GHI_{t}$ (global horizontal
irradiation) values are calculated using the deterministic cloud cover to $GHI$
model described in \cite{larson2016a} and given by \small
\begin{align}
  GHI_{t} = GHI_{t}^{CS} \cdot [0.35 + 0.65 \cdot (1-CC)]
\end{align}
\normalsize

\noindent where $cc$ is the normalized cloud coverage (0 = clear, 1 = overcast)
that is obtained from a \ac{NWP} model, and the $GHI_{t}^{CS}$ is the clear sky
$GHI$ that is estimated by the clear sky methods provided by the \emph{pvlib}
Python package \cite{Holmgren2018}.

\section{Case study}
\label{sec:case_study}
In this section, we describe the input data used by the \ac{HEMS} models
presented in Section \ref{sec:hems_setup}, and the technical details of the
simulation setup used to calculate the results.

\subsection{Electrical load}

\begin{figure*}
  \centering
  \includegraphics[width=\textwidth]{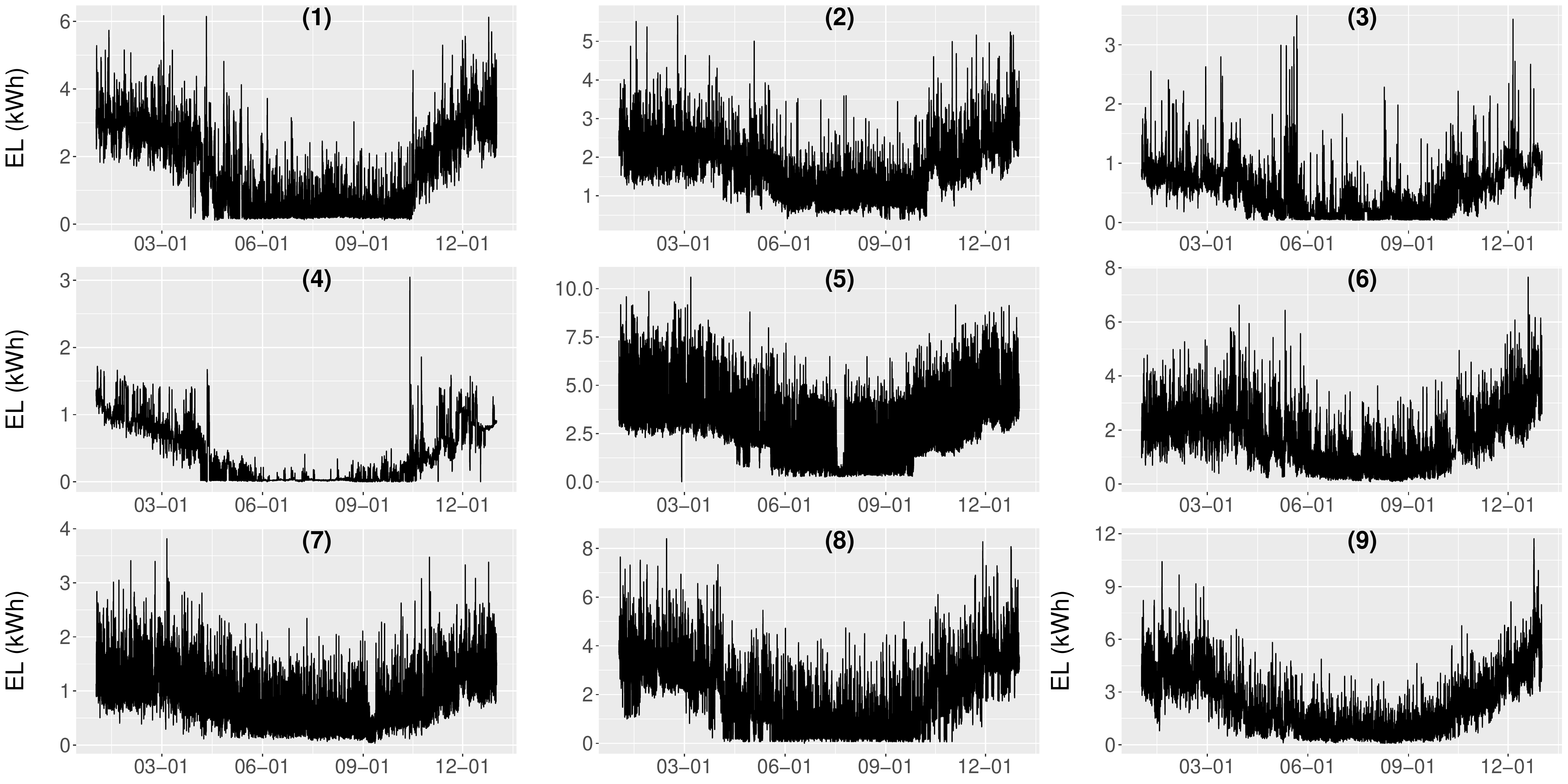} 
  \caption{Consumers' electricity consumption for the time period 2020-01-01 to
    2020-12-31.}
  \label{fig:users_ts}
\end{figure*}

The \ac{EL} demand profiles of nine residential consumers are shown in
Figure~\ref{fig:users_ts}. The consumption data results from smart meters
sampled at an hourly resolution for the year 2020. Information given about these
consumers includes the number of inhabitants, the approximate house location
given by its longitude and latitude coordinates, and the fact that they use heat
pumps as heating technology. The use of heat pumps explains the seasonality of
the \ac{EL}, i.e., a significantly higher consumption during winter in
comparison to summer (see Figure~\ref{fig:users_ts}). Although not visible in
yearly plots, intra-day patterns can be found on a closer inspection of the data
(see Figure~\ref{fig:user_data_set}). These patterns may be explained by the
daily routine activities of the tenants, e.g., having breakfast and dinner at
regular times. These factors together with the data presented in Section
\ref{sec:weahter_data} were considered when building the \ac{PLF} used for the
scenario generation of each consumer load demand. While it is out of the scope
of this paper to describe and analyze the \ac{PLF} methods, they are described
in detail in \cite{Lemos-Vinasco2021a}.

\subsection{Numerical weather prediction}
\label{sec:weahter_data}

Both the \ac{PV} simulation model (Section \ref{sec:pv_model}) and the \ac{PLF}
methods (Section \ref{sec:load_forecast_description}) rely on the \ac{NWP}
values as their primary input parameters. Here, the weather forecast provided by
the OpenWeatherMap service at an hourly resolution was used. The \ac{NWP} data
are described in \cite{OpenWeatherMap}. In particular, the ambient temperature
and solar irradiation were used for the \ac{PLF} and \ac{PV} models. Please note
that the solar irradiation signal was derived using the expected cloud coverage
from the \acp{NWP} in combination with the Global Horizontal Irradiation (GHI),
as described in Section \ref{sec:pv_model}.

\begin{figure}[th]
  \centering
  \includegraphics[width=8.5cm]{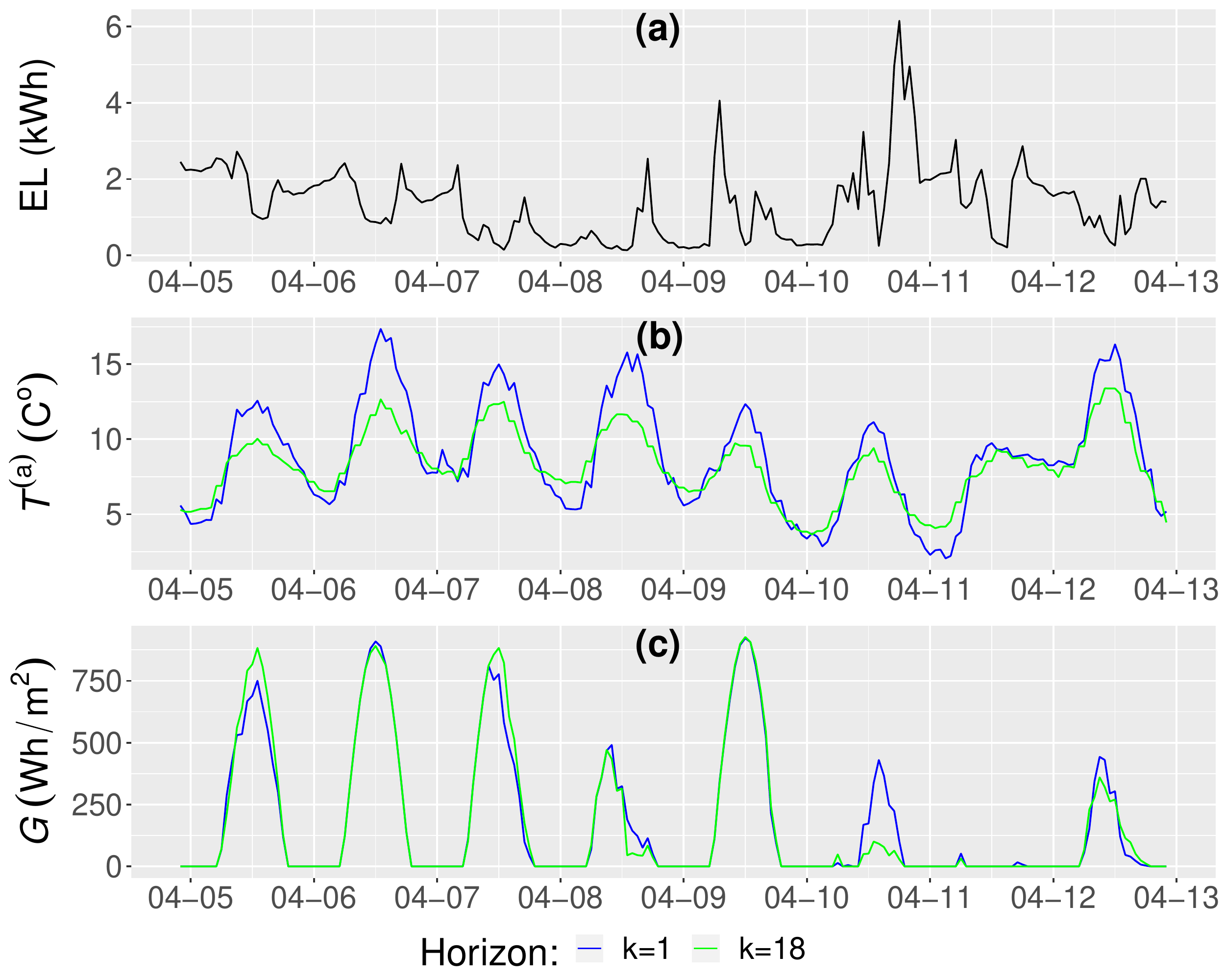}
  \caption{Consumer 1's electricity consumption (a), ambient temperature (b), and
    solar irradiation (c) for the time period 2020-04-05 to 2020-04-13.}
  \label{fig:user_data_set}
\end{figure}

The \ac{EL} and \ac{NWP} data are combined to produce a coherent dataset used
for the \ac{HEMS}. An example week for one user is shown in
Figure~\ref{fig:user_data_set}. Please note that the \ac{NWP} is updated every
hour with the forecast values covering several hours ahead ($k$ horizons). The
\ac{NWP} for one hour ahead $k=1$ and eighteen hours ahead $k=18$ are presented
in Figure~\ref{fig:user_data_set} (b) and (c).

\subsection{Electricity prices}

\begin{figure}[h]
  \centering
  \includegraphics[width=8.5cm]{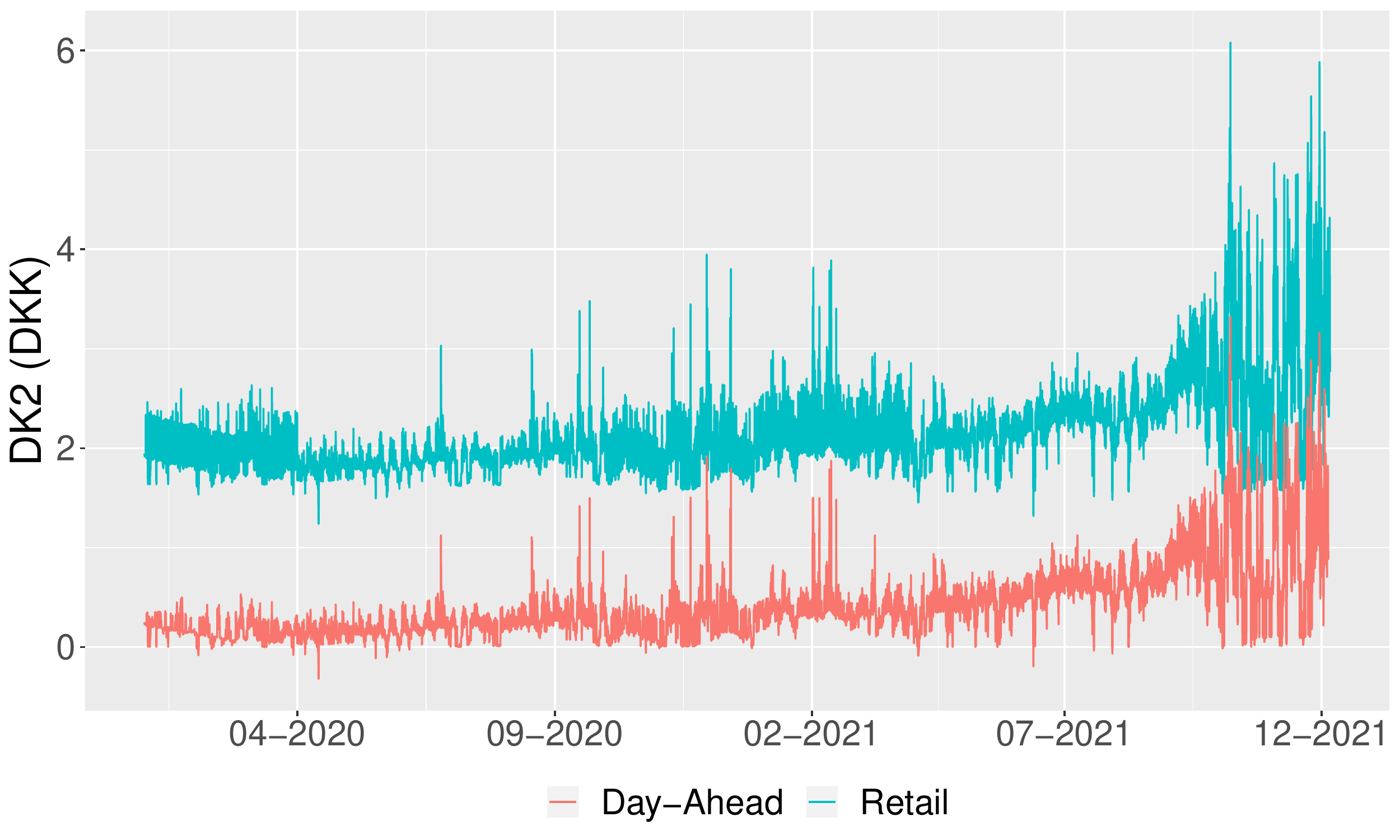}
  \caption{Retailer prices and \ac{DA} electricity prices for the time period
    2020-01-01 to 2021-12-03 for DK2.}
  \label{fig:prices}
\end{figure}

Given the high penetration of smart meters in Denmark, it is common for
electricity retailers to offer hourly prices to residential consumers.
Typically, this type of tariff is derived from the \ac{DA} wholesale electricity
market prices, also known as ELSPOT \cite{NordPool_DA}. In the ELSPOT market,
different zones/regions have their unique \ac{DA} prices. In Denmark, two price
zones zones exist: Western Denmark (DK1) and Eastern Denmark (DK2)
\cite{NordPool_areas}. To obtain the electricity prices for the residential
consumers, retailers add taxes, levies, and fees to the DA prices. In this
study, actual retail electricity prices provided by the Danish electricity
retailer Watts are used \cite{Watts_prices}. All consumers are located in the
greater Copenhagen area, which is a part of DK2 region. Figure~\ref{fig:prices}
shows the retailer prices and \ac{DA} electricity prices for the period of
2020-01-01 to 2021-12-03. The current Danish regulations allow residential
consumers to sell their surplus electricity back to the grid. The feed-in-tariff
is decided by the retailers. Most of them offer the ELSPOT price adjusted for
associated operational fees as feed-in-tariff to residential consumers, as
described by \cite{Vindstod}.

In our case study, Watts electricity prices are used as real-time price paid by
the consumers to purchase electricity, while \ac{DA} prices are used as
feed-in-tariff for selling. This corresponds to the parameters
$\ElecPrice{t}{+}$ and $\ElecPrice{t}{-}$, respectively (see Section
\ref{sec:hems_math}).

Please note the differences in prices in 2020 and 2021. On the one hand, a
tariff regulation change was introduced in 2021 that stipulates a low (between
00.00 to 17.00 and 20.00 to 00.00) and a peak (between 17.00 to 20.00)
electricity distribution fee from October to March \cite{Rad}. On the other
hand, 2021 was a year with unusual electricity prices, i.e., high prices and
high volatility \cite{Europeancomision}, as can be seen in the final quarter of
2021 in Figure \ref{fig:prices}.

The \ac{HEMS} setup and formulation allows cost reduction by selling excess
electricity (excess electricity from \ac{PV}), trading (buy at low-price hours
to sell at high-price ones), and load shifting. This can be done by exploiting
price volatility. Thus, for further reference, the mean and standard deviation
of the prices are presented in Table \ref{tab:price_volatility}. Note that the
statistics are only presented for January, April, July, and October of 2021.
These are the months that are included in the simulation setup described in
Section \ref{sec:simulation_study}. Please note that at present, consumers are
mainly passive users of electricity, which has no significant effect on prices
as they are price-takers \cite{Immonen2020}.

\begin{table}[h]
  \centering
  \footnotesize
  \begin{tabular}{c|rr|rr}
    \toprule
    \rowcolor{GrayDark}
    Time & \multicolumn{2}{c|}{Purchase price} & \multicolumn{2}{c}{Selling price}\\
    \rowcolor{GrayDark}
    Period & \multicolumn{1}{c}{Mean} & \multicolumn{1}{c|}{SD} & \multicolumn{1}{c}{Mean} & \multicolumn{1}{c}{SD} \\
    \hline
    \rowcolor{GrayLight}
    January 2020 & 1.943 & 0.164 & 0.206 & 0.074 \\
    \rowcolor{GrayLight}
    January 2021 & 2.085 & 0.230 & 0.379 & 0.126 \\
    April 2020 & 1.797 & 0.107 & 0.129 & 0.086 \\
    April 2021 & 2.008 & 0.120 & 0.357 & 0.160 \\
    \rowcolor{GrayLight}
    July 2020 & 1.858 & 0.144 & 0.191 & 0.115 \\
    \rowcolor{GrayLight}
    July 2021 & 2.330 & 0.198 & 0.615 & 0.159 \\
    October 2020 & 1.921 & 0.193 & 0.202 & 0.118 \\
    October 2021 & 2.624 & 0.780 & 0.810 & 0.591 \\
    \hline
  \end{tabular}
  \normalsize
  \caption{Purchase and selling prices mean and standard deviation (SD) for
    January, April, July and October.}
  \label{tab:price_volatility}
\end{table}




\subsection{Simulation setup}
\label{sec:simulation_study}

\begin{figure}[h]
  \centering
  \resizebox{8.5cm}{!}
  {
    \begin{tikzpicture}
      \tikzstyle{every node}=[font=\small]
      \draw[->] (0,0) -- coordinate (x axis mid) (11,0);

      \foreach \x in {0,...,10}
      \draw (\x,1pt) -- (\x,-4pt) node[anchor=north] {};

      \draw (1,1pt) (1,-4pt) node[anchor=north] {$t-N$};
      \draw (3,1pt) (3,-6pt) node[anchor=north] {$\dots$};
      \draw (5,1pt) (5,-4pt) node[anchor=north] {$t-2$};
      \draw (6,1pt) (6,-4pt) node[anchor=north] {$t-1$};
      \draw (7,1pt) (7,-4pt) node[anchor=north,color=blue] {$t$};
      \draw (8,1pt) (8,-4pt) node[anchor=north,color=blue] {$t+1$};
      \draw (9,1pt) (9,-6pt) node[anchor=north] {$\dots$};
      \draw (10,1pt) (10,-4pt) node[anchor=north] {$t+24$};

      \draw (4,4pt) (4,4pt) node[anchor=south]
      {Electrical load};

      \draw (4,20pt) (4,20pt) node[anchor=south]
      {Numerical weather predictions};

      \draw (8.5,4pt) (8.5,4pt) node[anchor=south]
      {Electricity prices};

      \draw (8.5,20pt) (8.5,20pt) node[anchor=south]
      {\ac{PLF}, \ac{PV}, and };

      \draw (7.8,-18pt) (7.8,-18pt) node[anchor=north,color=blue]
      {Operational};

      \draw (7.8,-28pt) (7.8,-28pt) node[anchor=north,color=blue]
      {set points};

      \draw [decorate,decoration={brace, mirror, amplitude=10pt, raise=1pt}]
      (1,-0.5) -- (7,-0.5) node [midway, anchor=north, outer sep=10pt]
      {Sliding window for forecast};

      \draw [decorate,decoration={brace, amplitude=10pt, raise=1pt}]
      (1,1.2) -- (7,1.2) node [midway, anchor=south, outer sep=10pt]
      {Values for \ac{PLF} estimation};

      \draw [decorate,decoration={brace, amplitude=10pt, raise=1pt}]
      (7,1.2) -- (10,1.2) node [midway, anchor=south, outer sep=10pt]
      {Rolling horizon for optimization};
    \end{tikzpicture}
  }
  \caption{Graphical representation of the simulation setting for time $t$.}
  \label{fig:sliding_window}
\end{figure}
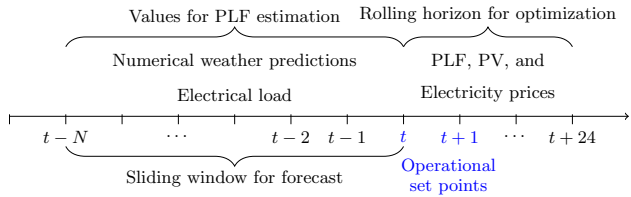

The simulation study is designed to resemble a real-time application. The aim of
the simulation is to optimize the battery's operational setpoints for the next
hour when considering a 24-hour horizon. Therefore, a rolling horizon approach
is used, which means that the \ac{PLF}, \ac{PV} simulation, and \ac{HEMS}
optimization will be updated every hour to determine the new operation
schedules. A graphical representation of the rolling horizon simulation setting
at time $t$ is presented in Figure \ref{fig:sliding_window}. Please note the
\ac{PLF} models are re-fitted using $t-N$ historical values at each time step
$t$. Thus, more accurate prediction can be expected by using the latest
available information from the forecasting models.

Four months of data (January, April, July and October) are selected as
representative for seasonal variations in order to analyse one year of
operation. Moreover, the prices used in the simulation are the prices in DK2
from 2021. This is motivated by the regulation changes introduced in that year.
Please note that the consumers' \ac{EL} data were only available for 2020.
Therefore, we used consumers' data from 2020 with prices from 2021 in our
simulations, assuming that the \ac{EL} in the selected months of 2020 is likely
to be similar to the \ac{EL} in the same months in 2021 and the fact that
residential consumers are price-takers.

\section{Simulation results}
\label{sec:results}
In this section, we compare consumers' cost savings when using different control
strategies. Two such strategies are considered: a naive controller and an
optimization-based controller. A naive controller refers to a consumer with
\ac{PV} and battery system without a \ac{HEMS}. This controller maximizes self
consumption by only selling electricity to the grid when the battery is fully
charged. The naive controller uses neither forecasting nor optimization methods
and it is usually the default controller in the \ac{HEMS} setup presented in
Section \ref{sec:hems_setup} \cite{Azuatalam2019a}. The optimization-based
controller refers to a consumer using a \ac{HEMS} with optimization and forecast
capabilities as presented in this paper. With this in mind, the rest of the
results section is organized as follows. In
Subsection~\ref{sec:methods_comparison}, we determine the most suitable \ac{PLF}
method for the optimization-based controller by comparing the performance of the
proposed \ac{HEMS} optimization model using different forecasting methods. A
comparison between the best optimization-based controller and a naive controller
is presented in Subsection \ref{sec:cost_results}. Results indicate the optimal
strategy is a combination of a naive and a optimal controller. This is presented
in detail in Subsection \ref{sec:naive_and_optimization_res}.


\subsection{Comparison of different optimization-based controllers}
\label{sec:methods_comparison}

\begin{table*}[t!]
\centering
\footnotesize
\scalebox{0.85}{
\begin{tabular}{c|rr|rr|rr|rr|rr}
\toprule
\rowcolor{GrayDark}
  Con. & \multicolumn{2}{c|}{PI-RH} & \multicolumn{2}{c|}{Copula-EV} & \multicolumn{2}{c|}{Copula-SP} & \multicolumn{2}{c|}{RLS-EV} & \multicolumn{2}{c}{RLS-SP} \\
  \rowcolor{GrayDark}
  no. & \multicolumn{1}{c}{DKK} & \multicolumn{1}{c|}{\%} & \multicolumn{1}{c}{DKK} & \multicolumn{1}{c|}{\%} & \multicolumn{1}{c}{DKK} & \multicolumn{1}{c|}{\%} & \multicolumn{1}{c}{DKK} & \multicolumn{1}{c|}{\%} & \multicolumn{1}{c}{DKK} & \multicolumn{1}{c}{\%} \\
  \midrule
  1 & 5798.61 & - & 6220.66 & 7.28 & 5976.54 & 3.07 & 6143.84 & 5.95 & 5895.99 & 1.68 \\ \rowcolor{GrayLight}
  2 & 12142.80 & - & 12887.78 & 6.14 & 12553.15 & 3.38 & 12745.82 & 4.97 & 12438.41 & 2.43 \\
  3 & 6689.17 & - & 7259.67 & 8.53 & 6977.89 & 4.32 & 7398.95 & 10.61 & 6984.19 & 4.41 \\\rowcolor{GrayLight}
  4 & 2071.31 & - & 2438.35 & 17.72 & 2210.56 & 6.72 & 2514.62 & 21.40 & 2183.16 & 5.40 \\
  5 & 1247.45 & - & 1638.04 & 31.31 & 1445.67 & 15.89 & 1695.82 & 35.94 & 1402.20 & 12.41 \\\rowcolor{GrayLight}
  6 & 7435.47 & - & 8322.06 & 11.92 & 7887.49 & 6.08 & 8254.65 & 11.02 & 7824.03 & 5.23 \\
  7 & 9294.46 & - & 9945.68 & 7.01 & 9610.14 & 3.40 & 9903.44 & 6.55 & 9557.67 & 2.83 \\\rowcolor{GrayLight}
  8 & 3157.25 & - & 3547.44 & 12.36 & 3312.33 & 4.91 & 3639.82 & 15.28 & 3316.83 & 5.05 \\
  9 & 6616.65 & - & 6970.68 & 5.35 & 6752.09 & 2.05 & 6954.68 & 5.11 & 6715.13 & 1.49 \\\rowcolor{GrayLight}
  \bottomrule
\end{tabular}
}
\normalsize
\caption{Total cost of the optimization and forecasting methods. All percentages
  are calculated relative to the PI-RH method.}
\label{tab:method_cost_comparison}
\end{table*}

The \ac{PLF} methods presented in \cite{Lemos-Vinasco2021a} allow different
combinations of forecasting and optimization methods. This section discusses the
performance of such combinations in order to select the best method. The
analyzed combinations are the following:

\begin{itemize}
  \item \textbf{RLS-SP}: the proposed \ac{HEMS} optimization using 100 scenarios
        generated by the RLS forecasting method.
  \item \textbf{RLS-EV}: the proposed \ac{HEMS} optimization using the expected
        value of the 100 scenarios generated by the RLS forecasting method.
  \item \textbf{Copula-SP}: the proposed \ac{HEMS} optimization using 100
        scenarios generated by the Copula forecasting method.
  \item \textbf{Copula-EV}: the proposed \ac{HEMS} optimization using the
        expected value of the 100 scenarios made by the Copula forecasting
        method.
  \item \textbf{PI-RH}: perfect information (PI) in a rolling horizon, i.e,
        using the proposed \ac{HEMS} optimization assuming that the consumer's
        load is known. This method is not applicable in practice, since the PI-RH
        method assumes a perfect knowledge of the future demand. But it can be
        used to give performance bounds on the optimization in the other
        settings.
\end{itemize}

Please note that ``EV'' methods correspond to a deterministic version of the
\ac{HEMS} stochastic model presented in Section \ref{sec:hems_math}. This allows
us to assess the impact of modeling the \ac{EL} uncertainty by comparing with
``SP'' methods. Table \ref{tab:method_cost_comparison} presents the total
electricity cost for the simulated months for the different combinations of
forecasting and optimization methods. Moreover, the table reports the cost
relative to the theoretical approach PI-RH expressed as percentages. The
simulation results indicate that the optimization using scenarios-based
stochastic programming outperforms the other solutions for all consumers, with
the RLS-SP being the best method. It presents the smallest difference to the
PI-RH. Considering that we re-optimize the solution after one time period, these
results are aligned with the results in \cite{Lemos-Vinasco2021a}. In
\cite{Lemos-Vinasco2021a}, although the Copula-based forecast presented an
overall better performance, the RLS-based forecast showed higher accuracy for
the initial time period. This may indicate that the performance of the
optimization depends to a high degree on the precision of the forecast in the
initial time step.

\begin{figure}[h]
  \centering
  \includegraphics[width=\textwidth]{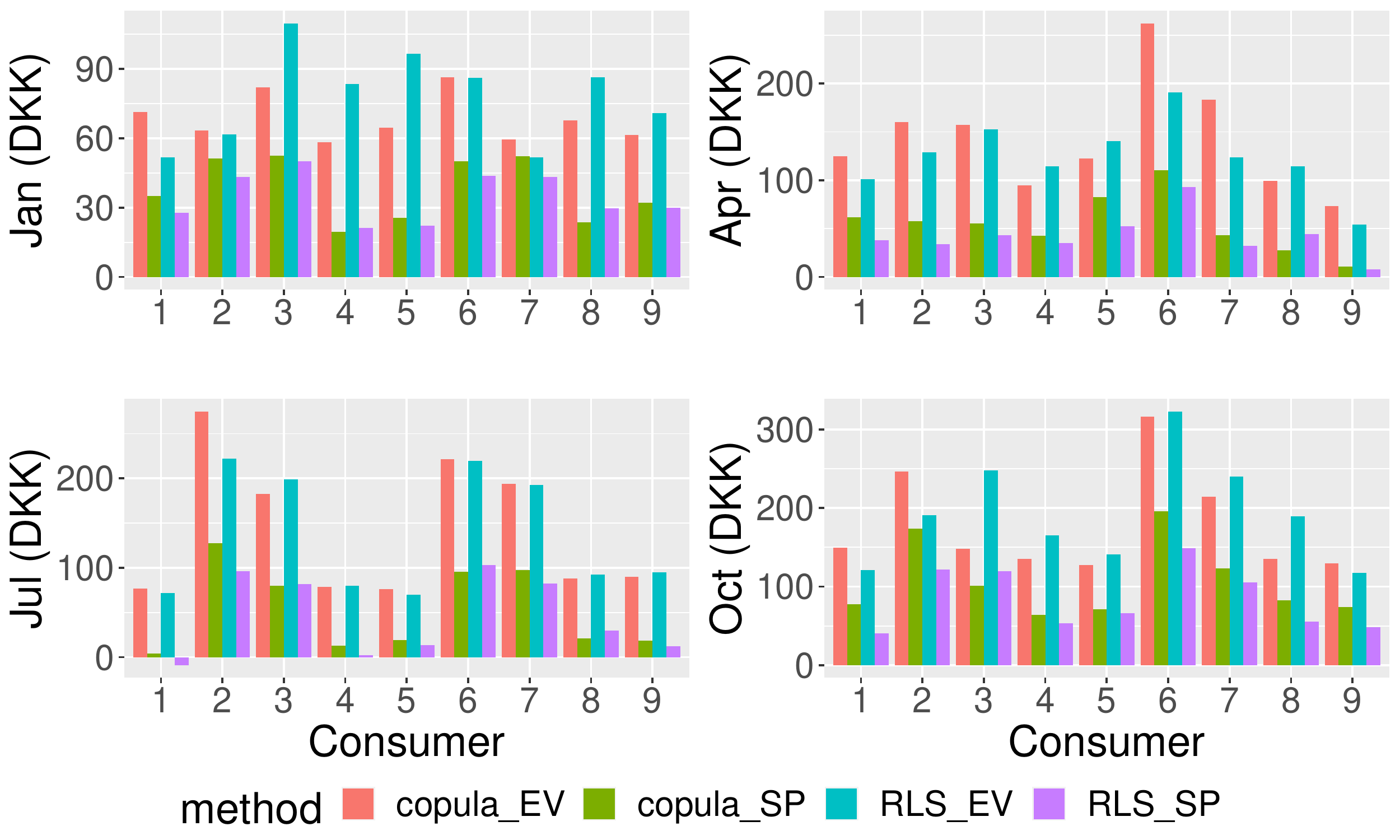}
  \caption{Additional cost by forecast and optimization framework in comparison
    with the PI-RH as baseline.}
  \label{fig:monthly_cost}
\end{figure}

A cost comparison between the different methods on a monthly level is shown in
Figure \ref{fig:monthly_cost}. The figure shows the additional cost incurred in
each method in comparison to the PI-RH case as our baseline. Similar results as
the aggregated results presented in Table \ref{tab:method_cost_comparison} are
seen, where the RLS-SP and Copula-SP methods outperform the deterministic
methods in every season. Moreover, the simulation results indicate that in some
cases, the solutions found through SP methods tend to be more robust than those
determined assuming perfect information. This can be seen in the RLS-SP results
for consumer 1 in July, where RLS-SP solutions outperform those of PI-RH. This
may seem counterintuitive considering that PI-RH assumes perfect knowledge of
uncertain \ac{EL}, however, the PI is limited by the 24 hours in the rolling
horizon. This behavior may indicate that considering the \ac{EL} stochasticity
may lead to more robust solutions towards the unpredicted \ac{EL} and/or the
optimization may benefit from a longer forecasting horizon.


\subsection{Comparison of naive and optimal control}
\label{sec:cost_results}

\begin{figure*}[ht]
  \centering
  \includegraphics[width=\textwidth]{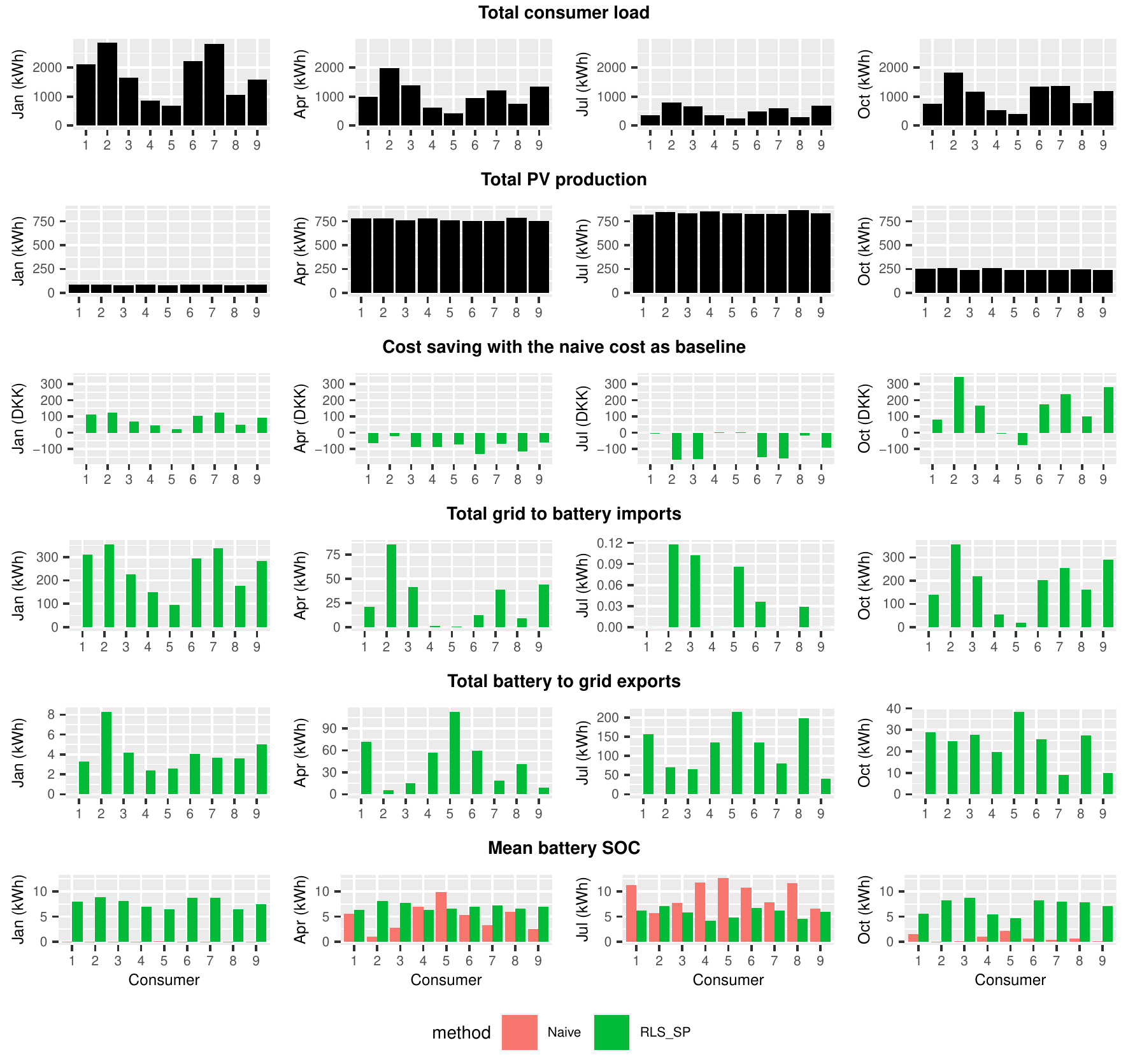}
  \caption{Monthly cost and operation results for the naive controller and RLS-SP.}
  \label{fig:operational_results}
\end{figure*}

In this section, we present the operational results for the optimization-based
control method (RLS-SP) and the naive controller. Figure
\ref{fig:operational_results} shows the total load, total \ac{PV} production,
cost savings relative to the naive controller, total amount of electricity
imported from the grid to the battery, the total power exported from the battery
to the grid, and the average battery SOC for each consumer.

From the cost saving plots, we can observe significant differences between
seasons. The RLS-SP outperforms the naive controller in winter (Jan) and autumn
(Oct). In January, the additional cost savings can be explained by more load
shifting. Here, the lower \ac{PV} production available for self-consumption
incentivizes more transactions with the grid, which can be seen by the higher
grid-to-battery energy flow. Note that the battery-to-grid energy flow is
non-significant. This indicates the optimization is taking advantage of the
price volatility (see Table \ref{tab:price_volatility}) by charging during
low-price hours in order to use the stored electricity during high-price hours.
October yields the highest cost savings of the simulated months. This is
achieved by exploiting the high price volatility (the highest price volatility
among the simulated months) in a similar fashion to January.

In contrast, spring (Apr) and summer (Jul) show the lowest cost savings. In
these months, the RLS-SP controller is outperformed by the naive strategy. In
particular, April presents a considerable \ac{EL} with high \ac{PV} production
but low price volatility (see Table \ref{tab:price_volatility}). This implies
that the optimization minimizes cost by selling excess \ac{PV} generation, which
could be a sub-optimal strategy in the long-term given the gap between purchase
and sale electricity prices (see Figure \ref{fig:prices}). In this case, the
daily rolling optimization horizon might not be able to capture the longer-term
effects and saving \ac{PV} generation for posterior use will benefit the
consumer the most. Furthermore, this argument matches the behavior seen in July,
where the low \ac{EL}, high \ac{PV} production, and low price volatility leave
almost no space for additional cost savings relative to the naive approach.
Therefore, under these conditions, a naive approach is able to minimize
consumers' costs in the long-term, and the application of more sophisticated
optimization-based methods could have a negative impact on cost minimization.


\subsection{Combination of naive controller and optimization}
\label{sec:naive_and_optimization_res}

As we have seen previously, the naive controller (which maximizes self
consumption) performs better than the optimization-based methods in the presence
of high \ac{PV} production and low \ac{EL}. This behavior could be explained as
a result of the price structure and volatility, and the implications of the
finite optimization horizon. The optimization minimizes short-term cost by
selling excess electricity, since it plans no more than 24 hours ahead. However,
sale prices are very unfavorable in comparison with purchase prices, which makes
self-consumption a better long-term strategy. Thus, one option could be to have
a \ac{HEMS} switching between a naive controller in the spring and summer and to
use an optimization-based control strategy in the autumn and winter seasons.
Hence, a full cost comparison between a passive consumer (without PV and
battery), a naive controller, and the proposed strategy (Naive+RLS-SP) is
presented in Table \ref{tab:total_cost_results}. The results show that
significant cost savings can be achieved in the naive and the proposed strategy
in comparison with a passive consumer. In particular, the simulation results
show that consumers with higher \ac{EL} benefit most from installing the
hardware and controllers. Moreover, the differences between the two control
strategies show that the combined controller (Naive+RLS-SP) provides on average
8.05\% additional savings for consumers with higher load (excluding consumers 4
and 5 with the lowest load of all consumers) in comparison with the naive
controller, as can be seen in detail in Table \ref{tab:total_cost_savings}.

\begin{table}[h]
  \centering
  \footnotesize
  \begin{tabular}{c|rr|rr|rr}
    \toprule
    \rowcolor{GrayDark}
    Consumer & \multicolumn{2}{c|}{Passive} & \multicolumn{2}{c|}{Naive} & \multicolumn{2}{c}{Naive+RLS-SP}\\
    \rowcolor{GrayDark}
    no. & \multicolumn{1}{c}{DKK} & \multicolumn{1}{c|}{\%} & \multicolumn{1}{c}{DKK} & \multicolumn{1}{c|}{\%} & \multicolumn{1}{c}{DKK} & \multicolumn{1}{c}{\%} \\
    \midrule
    1 & 9219.20 & - & 6015.12 & 34.75 & 5825.17& 36.81\\
    \rowcolor{GrayLight}
    2 & 16734.67 & - & 12715.96 & 24.01 & 12248.38& 26.81\\
    3 & 10799.02 & - & 6972.54 & 35.43 & 6735.15& 37.63\\
    \rowcolor{GrayLight}
    4 & 5353.58 & - & 2140.78 & 60.01 & 2094.73& 60.87\\
    5 & 3901.01 & - & 1281.02 & 67.16 & 1332.69& 65.84\\
    \rowcolor{GrayLight}
    6 & 11208.06 & - & 7827.33 & 30.16 & 7544.91& 32.68\\
    7 & 13370.53 & - & 9695.89 & 27.48 & 9333.01& 30.20\\
    \rowcolor{GrayLight}
    8 & 6455.34 & - & 3331.51 & 48.39 & 3184.64& 50.67\\
    9 & 10806.51 & - & 6935.86 & 35.82 & 6563.42& 39.26\\
    \bottomrule
  \end{tabular}
  \normalsize
  \caption{Total cost of the simulated months for passive consumers, and naive
    and optimal (naive+RLS-SP) control strategies. All percentage values are
    calculated relative to the passive consumers' cost.}
  \label{tab:total_cost_results}
\end{table}

\begin{table}[h]
  \centering
  \footnotesize
  \begin{tabular}{c|rrrr}
    \toprule
    \rowcolor{GrayDark}
    Consumer No. & \multicolumn{1}{c}{Naive} & \multicolumn{1}{c}{Naive+RLS-SP} & \multicolumn{1}{c}{Difference} & \multicolumn{1}{c}{\%}\\
    \midrule
    1 & 3204.09 & 3394.03 & 189.95 & 5.93\\
    \rowcolor{GrayLight}
    2 & 4018.70 & 4486.29 & 467.59 & 11.64\\
    3 & 3826.48 & 4063.87 & 237.39 & 6.20\\
    \rowcolor{GrayLight}
    4 & 3212.80 & 3258.85 & 46.05 & 1.43\\
    5 & 2619.99 & 2568.32 & -51.67 & -1.97\\
    \rowcolor{GrayLight}
    6 & 3380.73 & 3663.15 & 282.42 & 8.35\\
    7 & 3674.64 & 4037.51 & 362.87 & 9.88\\
    \rowcolor{GrayLight}
    8 & 3123.83 & 3270.70 & 146.88 & 4.70\\
    9 & 3870.65 & 4243.09 & 372.44 & 9.62\\
    \bottomrule
  \end{tabular}
  \normalsize
  \caption{Naive and optimal control strategies (naive+RLS-SP) total cost savings and
  their difference in value and percentage.}
  \label{tab:total_cost_savings}
\end{table}

\section{Discussion}
\label{sec:discussion}
The simulation results from Section \ref{sec:methods_comparison} show that
\ac{EL} uncertainty has a significant impact on the economic performance of the
optimization-based strategies, with the stochastic methods outperforming the
deterministic methods. While it appears clear that including \ac{EL} uncertainty
increases the consumers' economic benefits, the lack of access to real \ac{PV}
data limits the ability to study the impact of modeling the uncertainty inherent
to the \ac{PV} generation on the system operation. To consider the \ac{PV}
uncertainty, one needs to account for the correlation between the \ac{PV} and
\ac{EL} time series, which could be technically challenging in a multivariate
setting. This could be a topic for future research.

A comparison between naive and optimal control strategies, presented in Section
\ref{sec:cost_results}, shows that the seasonal variations and the electricity
price structure have a significant impact on the performance of the control
strategies. During cold seasons with low \ac{PV} production and higher \ac{EL},
an optimization-based strategy is preferable in order to exploit price
volatility by load shifting, charging the battery at low-price hours and
discharging it during high-price hours. In warmer seasons with low price
volatility, lower \ac{EL} and \ac{PV} make the naive strategy a better option.
This could be explained by the fact that under these conditions, a long-term
self-consumption strategy might be more profitable, which the optimization-based
controller fails to capture due to its finite optimization horizon. Please note
that this dynamic is tied to the price structure, where high differences between
purchase and sale prices do not incentivize consumers to trade electricity. This
may have a direct impact on the consumers in studies where market entities such
as aggregators are present. Under current market conditions, aggregators will
have to offer prices that surpass the retail prices. Otherwise, the consumers
will not participate in the aggregation programs. Thus, studying different
business models for aggregators under realistic price structures could be a line
of future research.

The simulation results presented in Section \ref{sec:naive_and_optimization_res}
indicate that the overall best control among the studied strategies is a
combination of a naive and optimization-based controllers. This is shown by the
proposed Naive+RLS strategy that outperforms the naive controller by 8.05\% on
average (excluding users 4 and 5). Please note that these extra savings result
from a software update on the default system controller, which does not incur
additional cost to the consumers. Additionally, the results obtained for
customers 4 and 5 suggest a self-consumption strategy benefits the consumer the
most when the \ac{EL} is low, making more elaborated control strategies
unnecessary in practice.

As smart grids continue to develop and price schemes evolve into a real-time
structure, the proposed setup and control strategy shows that there is an
incentive for consumers to adopt setups such as the one presented in this paper.
Moreover, if we consider that a higher energy demand and a higher penetration of
\ac{RES} is expected in the residential sector, we could argue that these types
of automation system might become a much needed tool for consumers to protect
themselves against higher price volatility. Furthermore, although the proposed
\ac{HEMS} omits elements such as direct control of smart appliances and electric
vehicles, we would expect their inclusion to bring higher savings for the
consumer, which is a further element that could be explored in future research.
While the economic results presented here show potential economical advantages
of \ac{HEMS}, the calculations were made for 2021 only. Future scenarios of
electricity prices could be included and analyzed in future research in this
area.

\section{Conclusion}
\label{sec:conclusion}
In this paper, a \ac{HEMS} setup tailored for residential consumers under
current electricity market rules in Denmark was presented. The \ac{HEMS} was
formulated as a \ac{MILP} using stochastic programming, where the \ac{EL} is
treated as an uncertain parameter, which is modelled by multivariate
probabilistic forecast models. Moreover, the \ac{HEMS} formulation allows
selling surplus electricity back to the grid depending on the price incentives
and the re-optimization of the solution in a rolling horizon fashion.

The simulation results suggest that the optimization-based controllers, which
consider the \ac{EL} uncertainty, outperform their counterparts that consider
only the expected value of the \ac{EL} for all users and seasons investigated.
Specifically, by comparing the RLS-SP controller with the naive one, we found
that the seasonality of the data and prices have a profound impact on the
cost-savings of all users. One possible explanation could be that the finite
optimization horizon and the gap between purchase and sale prices lead to
sub-optimal solutions of the optimization-based controller in spring and summer,
which are seasons with high \ac{PV} production and low \ac{EL}. In these
seasons, a naive controller performs better. On the other hand, in autumn and
winter with low \ac{PV} production and high \ac{EL}, optimization-based
controllers exploit price volatility by shifting load to minimize electricity
cost. Thus, the proposed combination of a naive and optimization-based strategy
was shown to reduce the annual consumers' electricity cost in comparison with
the naive method as a year-around controller. These savings are specifically
significant considering that the proposed change of strategies in different
seasons may come as a software update with small cost for the consumers.
Finally, the case study showed that the proposed setup and control strategy
offer a significant step towards a comprehensive assessment of the potential of
\acp{HEMS}, indicating that the setup and control strategy could be developed in
real-world applications of \acp{HEMS}.

\section*{Acknowledgements}
This work was partially funded by Innovation Fund Denmark through the
project identified with case no. 8053-00156B.

\raggedleft
\bibliographystyle{unsrt}
\bibliography{main}


\end{document}